\documentclass[3p,onecolumn,12pt,sort&compress]{elsarticle}

\usepackage{amssymb}
\usepackage{multirow}
\usepackage{subcaption}
\usepackage{amsmath}
\usepackage{nomencl}
\usepackage{lscape}
\usepackage{booktabs}
\usepackage{array}
\usepackage{xcolor}
\usepackage{float}
\usepackage{url}

\usepackage[version=3]{mhchem} 

\usepackage[symbol]{footmisc}

\usepackage{cleveref}
\crefname{chart}{Chart}{Charts}
\crefname{section}{Section}{Sections}
\crefname{figure}{Fig.}{Figs.}
\crefname{graph}{Graph}{Graphs}
\crefname{scheme}{Scheme}{Schemes}
\crefname{equation}{Eq.}{Eqs.}
\crefname{table}{Table}{Tables}
\crefname{chapter}{Chapter}{Chapters}
\crefname{appendix}{}{}
\newcommand{\data}[1]{\textcolor{black}{#1}}

\journal{European Journal of Operational Research}

\begin{document}

\begin{frontmatter}



\title{Funplex: A Modified Simplex Algorithm to Efficiently Explore Near-Optimal Spaces}


\author[inst1]{Christoph Funke\fnref{coauth}}
\author[inst1]{Linda Brodnicke\fnref{coauth}}
\author[inst2]{Francesco Lombardi}
\author[inst1]{Giovanni Sansavini\corref{cor1}}

\fntext[coauth]{Shared first co-authors}
\cortext[cor1]{Corresponding author: sansavig@ethz.ch}

\address[inst1]{Institute of Energy and Process Engineering, ETH Zurich, 8092 Zurich, Switzerland}
\address[inst2]{Faculty of Technology, Policy and Management, TU Delft, 2600 Delft, Netherlands}

\begin{abstract}
Modeling to generate alternatives (MGA) is an increasingly popular method in energy system optimization. MGA explores the near-optimal space, namely, system alternatives whose costs are within a certain fraction of the globally optimal cost. Real-world stakeholders may prefer these alternatives due to intangible factors. Nonetheless, widespread MGA adoption is hampered by its additional computational burden. Current MGA methods identify boundary points of the near-optimal space through repeated, independent optimization problems. Hundreds of model runs are usually required, and such individual runs are often inefficient because they repeat calculations or retrace previous trajectories. In this study, we transcend such limitations by introducing a novel algorithm called Funplex, which uses methods from multi-objective Simplex to optimize many MGA objectives with minimal computational redundancy. For a simple linear-programming energy hub case study, we show that Funplex is \data{five} times faster than existing methods and yields higher-quality near-optimal spaces. Furthermore, sensitivity analyses suggest that Funplex scales well with the number of investment variables, making it promising for capacity planning models. The current proof-of-concept implementation based on a full multi-objective tableau may face memory and stability limitations for large models. Nonetheless, future developments based on more advanced versions of Simplex may overcome such barriers, thereby making MGA more accessible and standard among modeling teams.
\end{abstract}




\end{frontmatter}



\section{Introduction}
\label{sec:intro}

Redesigning energy systems to meet climate mitigation goals is urgent \cite{IEA2023Breakthrough2023}, yet planning decisions on this matter are complex and slow-paced. There are many viable clean energy technology options and countless ways to combine them into feasible system designs. Individual designs can differ in terms of reliability, cost, investment needs, and social factors, further complicating planning decisions.

Energy system optimization models are employed to unravel such complexity and accelerate planning \cite{Susser2021}. They allow modelers to investigate various system design options and identify the one with the lowest cost for society. In the past, such cost-optimal solutions were the cornerstones of most model analyses.  Nonetheless, there is growing evidence that many feasible system designs exist, for marginally higher costs, that may be preferable due to a better alignment with the preferences of stakeholders and decision-makers \cite{trutnevyte2016,lombardi2020policy, Sasse2020Regional2035}. The energy transition involves many heterogeneous actors, from citizens to system operators and legislative entities, each with conflicting preferences regarding the choice of technologies and their spatial distribution  \cite{McGookin2024AdvancingModelling}. A single, cost-optimal solution fails to accommodate such preferences, thereby potentially resulting in poorly accepted system designs \cite{lombardi2020policy, Voll2015TheSynthesis}.

Approaches known as “Modeling to Generate Alternatives”  (MGA) have been gaining traction in the field of energy systems design to surmount the limitations of conventional cost-optimization \cite{Pickering2022DiversitySystem, patankar2023land, Sasse2023AVulnerabilities}. MGA identifies a large sample of feasible system designs that are economically comparable to the cost-optimal solution. Stakeholders may then appraise non-monetary trade-offs based on their preferences and knowledge, and eventually come to a consensus on preferable designs. MGA was initially proposed by \citeauthor{brill1982modeling} \cite{brill1982modeling}, with application to water management problems. In the past decade, the method was applied to energy system planning models by \citeauthor{decarolis2011using} \cite{decarolis2011using} and \citeauthor{trutnevyte2016} \cite{trutnevyte2016}, where it has met unprecedented success.

 In its simplest formulation, MGA explores the so-called near-optimal option space by repeatedly solving the original problem with different objective functions and a cost constraint. Each objective function seeks solutions which are different from previously explored options, while the cost constraint enforces that the cost remains within a marginal relaxation of the minimum feasible cost. Beyond this basic formulation, many variants of MGA exist. New MGA methods tailor to the needs of next-generation energy system models with high spatial and temporal resolution. For such models, hundreds to thousands of relevant, economically comparable design options exist that conventional MGA methods cannot efficiently explore. Each new method has different strengths and weaknesses. Some variants generate a small number of maximally different options, while others exhaustively find all possible design options. 

The MGA methods differ in which objective functions and constraints they use to find system alternatives. Among other things, the methods (1) minimize penalties for previously explored alternatives \cite{DeCarolis2016ModellingModel}; (2) maximize some distance metric \cite{Price2017ModellingModels}; (3) randomize objective functions \cite{trutnevyte2016}; or (4) constrain solutions to contain arbitrary amounts of the desired variables \cite{neumann2023broad}. For instance, \citeauthor{lombardi2020policy} \cite{lombardi2020policy} propose a spatially explicit version of MGA, named SPORES. SPORES enables the efficient identification of different spatial configurations of technology deployment in addition to different technology mixes. It leverages parallel computations starting from different extremes of the near-optimal space.  SPORES has been applied to European-scale models and carbon-neutral configurations of the Italian power system \cite{lombardi2023redundant, Pickering2022DiversitySystem}. Then, \citeauthor{neumann2023broad} \cite{neumann2023broad} explore the near-optimal space by first determining a few variables of interest and then using an epsilon-constraint method to delineate the Pareto front between these factors. The idea is to obtain, in such a way, a more uniform sample of the option space. They also apply their method to a European model.  
Next, several versions of MGA exist which use randomly generated coefficients in the objective function   \cite{trutnevyte2016, Berntsen2017EnsuringAlternatives, Patankar2023LandWest}. These methods have remained popular, even on large models, because they provide good samples of the near-optimal space if a large number of runs can be carried out. Finally, \citeauthor{pedersen2021modeling} 
\cite{pedersen2021modeling}, describe a method of finding the complete convex hull of the near optimal space in a few chosen dimensions. After identifying the convex hull, they rely on sampling techniques to identify interior points within the region.  

Despite its merits and the recent methodological developments,  MGA is not widely adopted due to its high computational burden \cite{lombardi2023redundant}. Methods which compute the complete convex hull, such as \citeauthor{pedersen2021modeling} \cite{pedersen2021modeling}, are computationally infeasible for applications with more than ten variables of interest \cite{pedersen2021modeling}. This holds true even on powerful high-performance computing machines. Next, algorithms which seek maximally different solutions require long computational times to reproduce a good representation of the actual near-optimal space. Many recent algorithms of this kind allow for a partial parallelization of the search by leveraging different, independent starting points. Nonetheless, the computation remains time and resource-intensive for large models \cite{lombardi2023redundant}. Finally, algorithms based on random coefficients \cite{trutnevyte2016, Berntsen2017EnsuringAlternatives, Patankar2023LandWest} can be completely parallelized; however, in practice computational resources limit the amount of simulations that can be run simultaneously. 

Improving the computational efficiency of MGA is a high priority, and there is untapped potential to do so. Existing MGA methods solve repeated optimization problems to identify system alternatives. During each optimization, solvers communicate only the optimal solution to a given problem. Intermediary solutions are ignored, even though these are often valid near-optimal alternatives themselves. To obtain these alternatives, they would need to be later re-identified in a separate MGA run. Reprogramming solvers so that they store, utilize, and communicate the intermediary near-optimal solutions may thus bring computational benefits. 

With this study, we set out to investigate this untapped potential by devising a first-of-its-kind MGA-tailored solver algorithm.  In particular, the algorithm we propose, named Funplex, is an MGA-tailored advancement of the Simplex algorithm. It represents the first attempt at improving MGA’s computational feasibility by acting directly on solvers’ algorithms rather than on the mathematical formulation of the MGA problem only. The idea of reformulating solver algorithms in this direction has been proposed so far only for the specific class of mixed-integer problems \cite{Gurobi2024SolutionOptimization, IBM2024SolutionDocumentation, Danna2007GeneratingProblems}. These are not applicable to large, linear, state-of-the-art energy system models \cite{Danna2007GeneratingProblems}. 

Funplex resembles and builds upon developments in multi-objective Simplex (MS) which were piloted in the 70s \cite{Evans1973APrograms, Ecker1978FindingPrograms, Ehrgott2005MulticriteriaOptimization}. As underscored by \citeauthor{neumann2023broad}, the MGA problem is closely intertwined with multi-objective optimization because the boundary of the near-optimal space forms a Pareto front \cite{neumann2023broad}. MS attempts to find all vertices on exactly such a Pareto front. The algorithm starts by optimizing one objective to find a single Pareto-optimal vertex. Then, additional sub-optimization problems are solved to identify which edges to follow to remain on the Pareto front \cite{Ehrgott2005MulticriteriaOptimization}. This process is repeated, until all Pareto-optimal vertices are found. While MS determines the Pareto with high accuracy, its many sub-optimizations, even for small problems, are computationally intensive \cite{Ehrgott2005MulticriteriaOptimization}. Funplex forgoes this limitation by trading accuracy for computational efficiency. The algorithm does not seek to find all efficient vertices but simply tries to stay near the boundary.                                                           
We test Funplex's performance on a simple energy hub model and analyze its computational efficiency as well as the quality of near-optimal space it finds. During these tests, we benchmark Funplex against the above-mentioned SPORES \cite{lombardi2023redundant} and Random Directions \cite{Berntsen2017EnsuringAlternatives} algorithms. SPORES is among the most recent approaches for efficient computation of the near-optimal space, and Random Directions is very popular and assumed to ensure good sampling of the near-optimal space (for more details on each Algorithm, see \cref{sec:appendix_MGA_methods}). Our results show that Funplex is substantially more efficient than existing methods while also yielding a complete representation of the near-optimal space. Through a sensitivity analysis, we also show that Funplex scales well with the number of investment variables, making it promising for larger capacity planning models. At the same time, we detect memory and stability limitations that may counterpoise its benefits with application to large models. This is in line with the expectations since the current implementation of Funplex is meant as a proof-of-concept. Reflecting on the results, we elaborate on future developments to overcome such barriers, thereby paving the way for making MGA a more accessible and standardized modeling practice.

\section{Methods}
\label{sec:Methods}
In this study, we propose an algorithm called Funplex, which is a modified version of the primal Simplex algorithm. Funplex is described in \cref{sec:algorithm}, and the interested reader is referred to \cref{sec:appendix_Simplex} for a brief introduction to the primal Simplex algorithm. Funplex's performance is benchmarked against two established MGA methods, namely Random Directions and SPORES, which are considered state-of-the-art methods. Random Directions enables the most uniform sampling and is thus expected to perform well in near-optimal solution space quality. SPORES is an example of an effort to improve the computational efficiency of identifying the near-optimal space and is thus expected to perform well in terms of computational efficiency. \cref{sec:appendix_MGA_methods} describes these two established methods. The criteria for evaluating the three algorithms' performance in terms of near-optimal solution space quality and computational efficiency are detailed in \cref{sec:KPI}. A brief description of the case study on which the analysis is performed is included in \cref{sec:case_study}.

\subsection{Proposed algorithm}
\label{sec:algorithm}
\begin{figure}[!h]
    \centering
    \includegraphics[width = 0.9\textwidth]{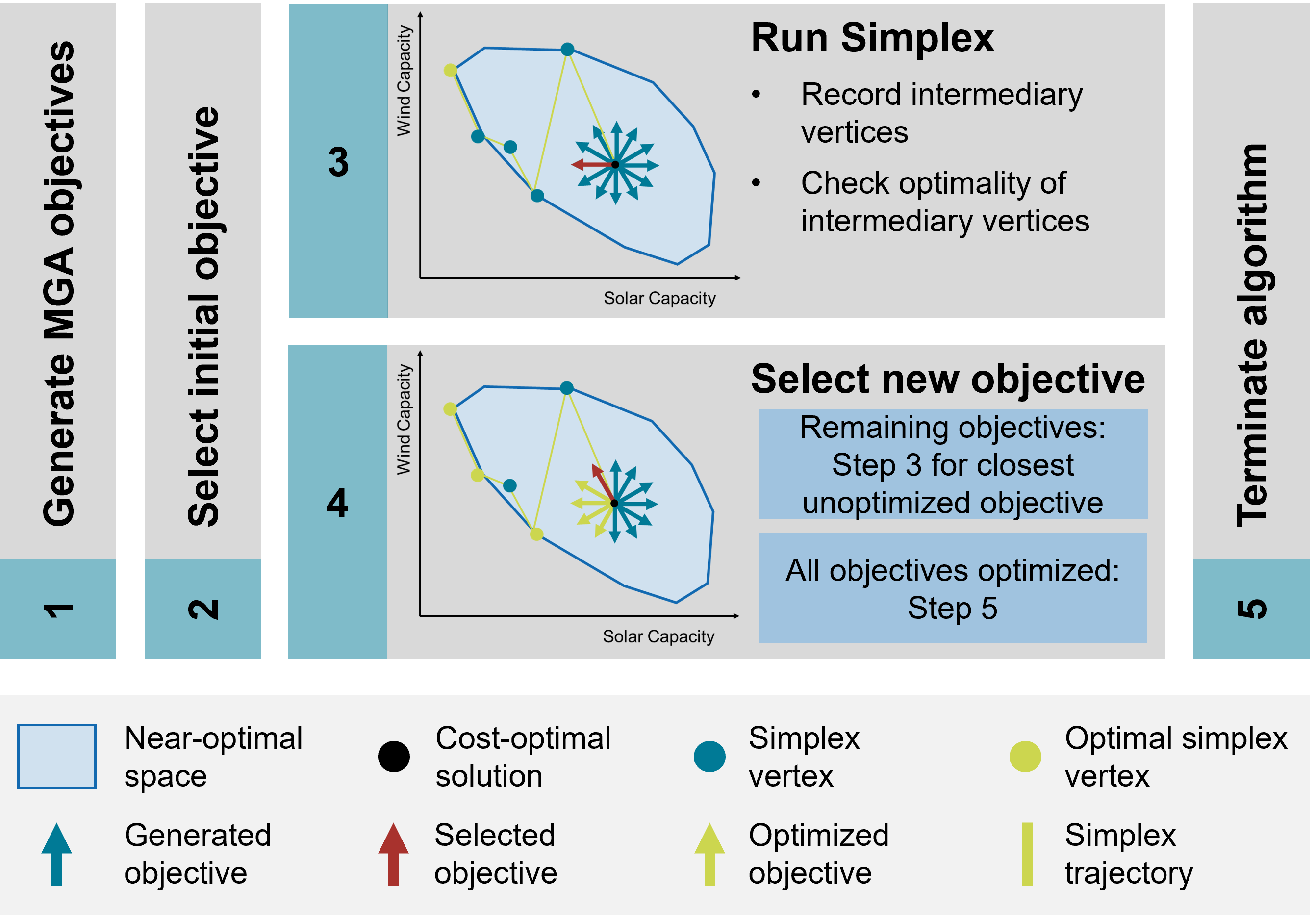}
\caption{Overview of Funplex's five steps illustrated using the two-dimensional illustration of the near-optimal space in terms of solar and wind capacity.}
    \label{fig:methodological_overview}
\end{figure}

The core idea of Funplex is to reduce computational redundancy while solving many objectives over the same feasible region. It does so by leveraging the fact that solvers, in this case, the primal Simplex algorithm, walk along the boundaries of the feasible space. Throughout each optimization, Funplex keeps track of intermediary solutions to gain additional information about the near-optimal space. Furthermore, it uses information from previously solved objectives to choose efficient starting points for optimizations. Funplex uses a multi-objective Simplex tableau to store and track all MGA objectives throughout this process efficiently. Such a tableau concurrently includes each objective's objective values and relative cost vectors. \cref{fig:methodological_overview} outlines the five steps of the Funplex algorithm, which we describe in depth in the following subsections.

In the following sections, let $N_d$ be the number of decision variables of interest $\mathbf{x^\ast} = [x_{1}^{\ast}, x_{2}^{\ast}, ..., x_{N_d}^{\ast}]$ for which we want to find a near-optimal space. In other words, $N_d$ is the dimension of the desired near-optimal space. Next, let $N_k$ be the desired number of objectives to explore the near-optimal space, i.e., the number of directions in which the near-optimal space is explored.

\subsubsection{Step 1: Generate MGA objectives}
Funplex begins by establishing a set of MGA objectives to optimize. For this step, any existing, non-iterative method of selecting MGA functions would work. Non-iterative methods pre-compute all MGA objectives prior to the analysis. Funplex generates these objectives from a random hypersphere. The following describes the procedure to obtain one objective function $f_{\mathbf{d_k}}$. To obtain the complete set of objectives $\mathcal{K}$, Funplex repeats this process $N_k$ times.

Points from a hypersphere can be computed from Gaussian random numbers. Let $\{z_1, z_2, ..., z_{N_d}\}$ be random numbers from a standard normal distribution. Then the point

\begin{equation}
    \mathbf{d} = 
        \begin{bmatrix}
           d_{1} \\
           d_{2} \\
           \vdots \\
           d_{N_d}
        \end{bmatrix} 
        = \frac{1}{\sqrt{z_1^2 + z_2 ^2 + ... + z_{N_d}^2}}
        \begin{bmatrix}
           z_{1} \\
           z_{2} \\
           \vdots \\
           z_{N_d}
        \end{bmatrix}
\end{equation}

\noindent is sampled randomly from the surface of an $N_d$ dimensional hypersphere \cite{Hicks1959, Marsalia1972,Muller1959}. Note that $\lVert\mathbf{d}\rVert = 1$ by definition. For $N_d = 2$, these points lay on a unit circle, whereas for $N_d = 3$, these points lay on a unit sphere.   

We can then define an MGA objective function for each direction vector ($\mathbf{d}$). The objective function can thus be expressed in terms of the decision variables of interest $\mathbf{x}^{\ast} = [x_{1}^{\ast}, x_{2}^{\ast}, ..., x_{N_d}^{\ast}]$:

\begin{equation}
    \label{eq:funplex_objective}
    f_{\mathbf{d}}= \sum_{i = 1}^{N_d} \frac{d_i}{L_i}x_{i}^{\ast}
\end{equation}

Here, $L_i$ is a characteristic scale for the decision variable of interest $x_i^{\ast}$. This characteristic scale approximately normalizes the decision variables, helping improve performance in cases where the decision variables span vastly different scales (e.g., 10 vs. 10,000). In this paper, the characteristic scales are taken from the values of $\mathbf{x}^{\ast}$ in the optimal solution, when available. When these are zero, the characteristic scales are estimated to roughly match the expected magnitude of the variables in the near-optimal space based on expert judgment. The characteristic scales are approximate, and they have not been tuned in this analysis. 

\subsubsection{Step 2: Select initial objective}
Funplex begins by selecting one objective arbitrarily. In this study, we select the first of the generated objectives.

\subsubsection{Step 3: Run Simplex}
Funplex uses the primal Simplex algorithm to optimize objectives with a few minor modifications. First, Funplex stores each intermediary solution, and second, Funplex evaluates the optimality of all objectives ($\mathcal{K}$) after each pivot. To do so, Funplex uses a multi-objective Simplex tableau. The Simplex implementation by Daneshian \cite{daneshian2023simplex} constitutes the groundwork from which we develop our Funplex.

In each pivot, Funplex determines the pivot column ($j$) by assessing the relative cost vector corresponding to the current objective $k \in \mathcal{K}$. There are multiple options for how to pick the pivot column from the relative cost vector \cite{zenklusen2023linear}. This analysis chooses the column that leads to the highest marginal decrease in the objective value, denoted by the most negative entry in the relative cost vector. This is standard for Simplex, although other methods exist \cite{zenklusen2023linear}. After the pivot column has been chosen, the pivot row is chosen among the legal options \cite{zenklusen2023linear}. Legal options lead to a new feasible vertex in which all decision variables are non-negative \cite{zenklusen2023linear}. Usually, there is only one legal pivot row given a chosen pivot column.

At each intermediary point, Funplex assesses the optimality of all $\mathcal{K}$ objectives by looking at their relative cost vectors in the tableau. If a relative cost vector has only zeros or positive elements, the corresponding objective has been optimized. Funplex tracks which objectives have been optimized using a boolean array $success \in \mathbb{R}^{N_d}$. A zero entry indicates that the corresponding objective has not been optimized, while a one indicates that it has been optimized.

\subsubsection{Step 4: Select new objective}
After the first objective has been optimized, Funplex selects a "similar" next objective. There are multiple similarity metrics, and this paper measures similarity using angle difference. Let $\mathbf{d_1}$ and $\mathbf{d_2}$ be the underlying directions for objectives $f_{\mathbf{d_1}}$ and $f_{\mathbf{d_2}}$, respectively. Then, we define the distance between the two objectives as:

\begin{equation}
    \label{eq:angle_difference}
    dist(f_{\mathbf{d_1}}, f_{\mathbf{d_2}}) = \mathbf{arccos}(\mathbf{d_1}\cdot\mathbf{d_2})
\end{equation}

\noindent This equation is derived from the definition of the angle ($\theta$) between two vectors:

\begin{equation}
    \label{eq:dot_product}
    \cos(\theta) = \frac{\mathbf{d_1}\cdot\mathbf{d_2}}{\lVert\mathbf{d_1}\rVert \lVert\mathbf{d_2}\rVert}
\end{equation}

\noindent Note that the directions $\mathbf{d_1}$ and $\mathbf{d_2}$ are unit vectors, therefore $\lVert\mathbf{d_1}\rVert = \lVert\mathbf{d_2}\rVert = 1$.

When choosing the next objective, Funplex searches the remaining (unoptimized) objectives for the one which is closest to the current objective by \cref{eq:angle_difference}. Since this objective is similar to the previous one, the difference between the two optima will likely be small. 

\subsubsection{Steps 5: Terminate algorithm}
The Funplex stops once each objective ($k \in \mathcal{K}$) has been optimized at least once. It then returns the coordinates of all optimal points and intermediary vertices.

\subsection{Performance evaluation}
\label{sec:KPI}
\subsubsection{Quality of near-optimal spaces}
MGA algorithms aim to explore the near-optimal space in all $N_d$ dimensions and discover as much of the "real" near-optimal space as possible. Thus, an ideal MGA algorithm discovers the largest possible near-optimal space with a reasonable computational effort. 

\paragraph{Two-dimensional projections} A qualitative comparison of the quality of the found spaces is enabled by assessing the two-dimensional projections of the near-optimal space. For the three algorithms, the boundaries of the identified near-optimal space are plotted and compared to the "real" near-optimal space which is computed via a two-dimensional planar MGA analysis, i.e., on only the two variables of interest, and with a very large number of objectives. We assume such a procedure to approximate the true near-optimal space with high fidelity. 

\paragraph{Volume of near-optimal spaces} A quantitative comparison of the three algorithms is enabled by calculating and comparing the volume of the $N_d$-dimensional space encapsulated by the convex hull of a given number of $N_k$ MGA points. To make the volumes easier to interpret, we normalize reported volumes so that the volume determined by Funplex has a size of 1. 

\paragraph{Volume gain} To enable a quantitative comparison of the three algorithms' performance in terms of solution quality, the volume gain is calculated. Volume gain measures the improvement resulting from Funplex for a given number of objectives $N_k$:
\begin{equation}
    \text{volume gain} = \frac{\text{Volume Funplex}}{\text{Volume established method}}
\end{equation}
Thus, a volume gain of 2 corresponds to doubling of the volume of the identified $N_d$-dimensional near-optimal space.

\subsubsection{Computational efficiency}
Three indicators are used to quantify the computational efficiency of each method, namely the number of tableau pivots, the computational complexity per pivot, and the runtime. Moreover, the term efficiency gain is defined, which is used to compare how the performance of Funplex scales with respect to the other two algorithms. All optimizations for SPORES and Random Directions are solved using a self-programmed MATLAB Simplex algorithm to enable a fair comparison since the value of each indicator quantifying computational efficiency depends on the exact solver implementation. The Simplex solver uses the same code as Funplex for finding feasible solutions, creating tableaus, and pivoting tableaus. 

\paragraph{Number of tableau pivots} The number of tableau pivots measures how many vertices each method transverses to find the solution. This is a useful performance indicator because the runtime of both Simplex and Funplex scales linearly with the number of pivots. Only pivots from the phase-two problem, after a feasible basic solution has been found, are counted.

\paragraph{Computational complexity per pivot} This performance indicator measures the number of calculations required per tableau pivot. We report the complexity in "Big O" asymptotic notation. It describes how the number of calculations changes as the size of the tableau increases, ignoring any constant factors. The number of calculations differs between Funplex and Simplex because Funplex uses a multi-objective tableau rather than a single-objective tableau.

\paragraph{Runtime} Runtime measures the time it takes for each method to run in our MATLAB code. The runtime includes finding an initial feasible solution, setting up a tableau, and optimizing the tableau. For SPORES and Random Directions, we reuse the same feasible solution between optimization problems but create a new tableau each time. This is consistent with what is currently feasible on commercial solvers. Generally, runtime depends on the exact algorithm implementation and used hardware. Although intuitive, it is, therefore, not a particularly robust indicator.  

\paragraph{Efficiency gain} The efficiency gain describes how the computational complexity of Funplex compares to that of SPORES and the Random Directions algorithm for a given number of objectives $N_k$. It is defined as:

\begin{equation}
    \text{efficiency gain} = \frac{\text{\# pivots established method}}{\text{\# of pivots Funplex}}
\end{equation}
Thus, an efficiency gain of ten corresponds to a tenfold reduction in required pivots.

\subsection{Case study}
\label{sec:case_study}

The case study is based on a linear model for designing an energy hub published by Hohmann \cite{hohmann2017energy}. The energy hub must supply heat and electricity to meet the baseload demand of a solid sorbent direct air capture (DAC) unit. The problem is solved for a representative operation year and results in \data{450} investment and operation decision variables and \data{494} constraints, as further detailed in the following paragraphs. It can be solved using Gurobi in \data{0.2} seconds and our MATLAB Simplex code in \data{7} seconds. To find the cost-optimal solution of the energy hub model, the overall system cost is minimized, which comprises of annualized investment and fuel import costs. 

\begin{figure}[h!]
    \centering
    \includegraphics[width = 0.8\textwidth]{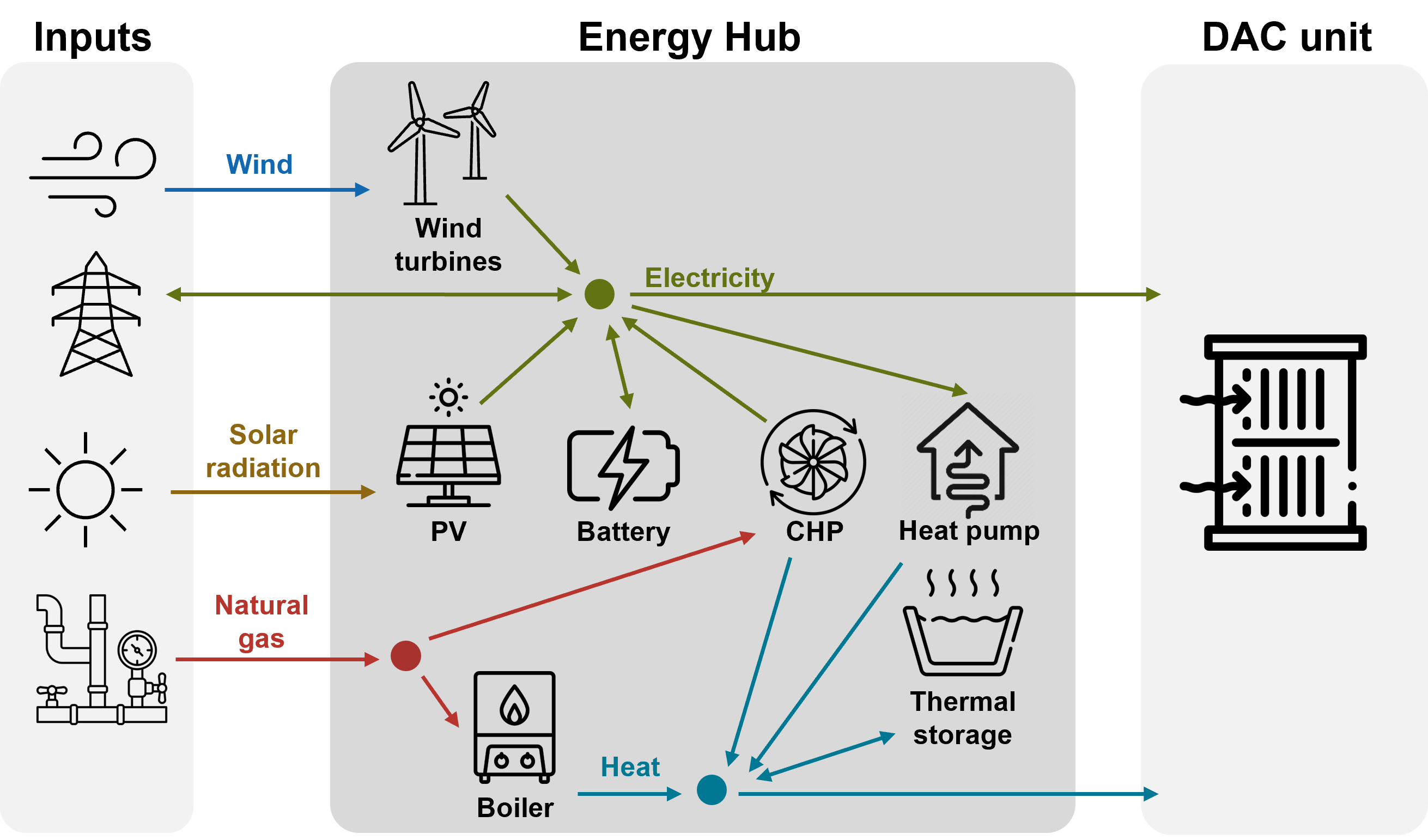}
\caption{Overview of the modelled energy system. The energy hub model optimizes the investment and operation of the seven available technologies over the time horizon of one year, which is approximated using a user-defined number of representative days.}
    \label{fig:MES_illustration}
\end{figure}

\paragraph{Technologies and input data}
The baseload heat and electricity are supplied at a single node through the installation and operation of seven available technologies, as well as via electricity and gas imports. The seven available generating technologies are wind turbines, photovoltaic (PV) panels, a battery, a combined heat and power (CHP) unit, a heat pump, a gas boiler, and a thermal storage (TS). The resulting reference energy system is illustrated in \cref{fig:MES_illustration}.

The investment costs for each technology are annualized, assuming a \data{25}-year lifetime. Electricity and gas imports are associated with energy import costs and a grid-average carbon emission rate. The heat demand is assumed constant at approximately \data{1.1 MW}, and the electricity demand is constant at approximately \data{440 kW}. Moreover, the system is subject to a carbon constraint, restricting yearly \ce{CO2} emissions to \data{1460} tons of \ce{CO2}. 

\paragraph{Time horizon} The considered time horizon for the operation is one year, which is modeled based on the number of representative days and the number of hours per day that the user selects. Representative days are generated by using a k-means clustering algorithm on the annual, hourly profile of solar irradiance, wind, heat demand, and electricity demand. The final profiles are taken from the cluster centroids and, therefore, represent average days throughout the year. One representative day with twenty-four hours is used for the base model. 

\paragraph{MGA assumptions}
The algorithms are tested for $N_d = 4$ decision variables of interest, namely the installed technology capacity of (1) wind, (2) the gas boiler, (3) PV, and (4) the heat pump. This number of dimensions is chosen since it allows for a meaningful illustrative MGA analysis whilst still being tractable to visualize. All analyses are performed for a slack value of 5 \%, i.e., a maximum cost deviation of 5 \% from the minimum-cost solution. 200 objective functions are used for all MGA algorithms. Since both Funplex and Random Directions rely on generating random numbers, the random number seed is fixed to ensure the reproducibility of the results. 

\subsection{Scalability analysis} Funplex is implemented as a proof-of-concept algorithm and tested on a small-scale problem. Thus, a scalability analysis is performed to evaluate how Funplex might perform on a larger model. The scalability is examined in two dimensions, namely the performance changes with (1) increasing model size and (2) increasing MGA complexity. 

\paragraph{Scalability with model size} 
Increasing model size is implemented by (1) increasing the time horizon of the model and (2) increasing the number of investment variables. The increase of investment variables is modeled by increasing the number of available solar sites, each of which can be invested in individually. This enables users to increase the number of investment variables and add a spatial component to the model. The solar irradiation profile for each solar site is taken by adding multiplicative log-normal noise to the base profile. The noise factor was taken from a log-normal distribution with a mean of 1 and a standard deviation of 0.254. When assessing the impact of model scalability with the number of PV sites, a reduced six-hour operating horizon is used to ensure numerical stability and allow a larger number of sites to be simulated.

\paragraph{Scalability with MGA complexity}
The scalability with MGA complexity is explored along two dimensions, i.e., (1) the number of MGA objectives and (2) the number of MGA dimensions. 
The number of MGA objectives is achieved by varying the number $N_k$. For Funplex and the Random Directions algorithm, this means changing the number of randomly generated objectives. For the SPORES algorithm, the number of explored objectives per spore is adjusted accordingly. The number of MGA dimensions is adjusted by varying the number of decision variables of interest $N_d$. This is done by including up to all seven buildable technologies in the model. 
\section{Results and discussion}
\label{sec:Results}
The following subsections present the results and discussion of our analysis, starting with a presentation of the cost-optimal model results in \cref{sec:cost_optimal}. \cref{sec:quality} provides a qualitative and quantitative comparison of the solution quality of the three algorithms, followed by a discussion of their performance in terms of computational efficiency in \cref{sec:efficiency}. \cref{sec:scalability} assesses the scalability of Funplex to larger models and more complex MGA analyses, and finally \cref{sec:limitations} discusses the limitations of our study.

\subsection{Cost-optimal solution}
\label{sec:cost_optimal}
\begin{figure}[ht]
    \centering
    \includegraphics[width = 0.5\linewidth]{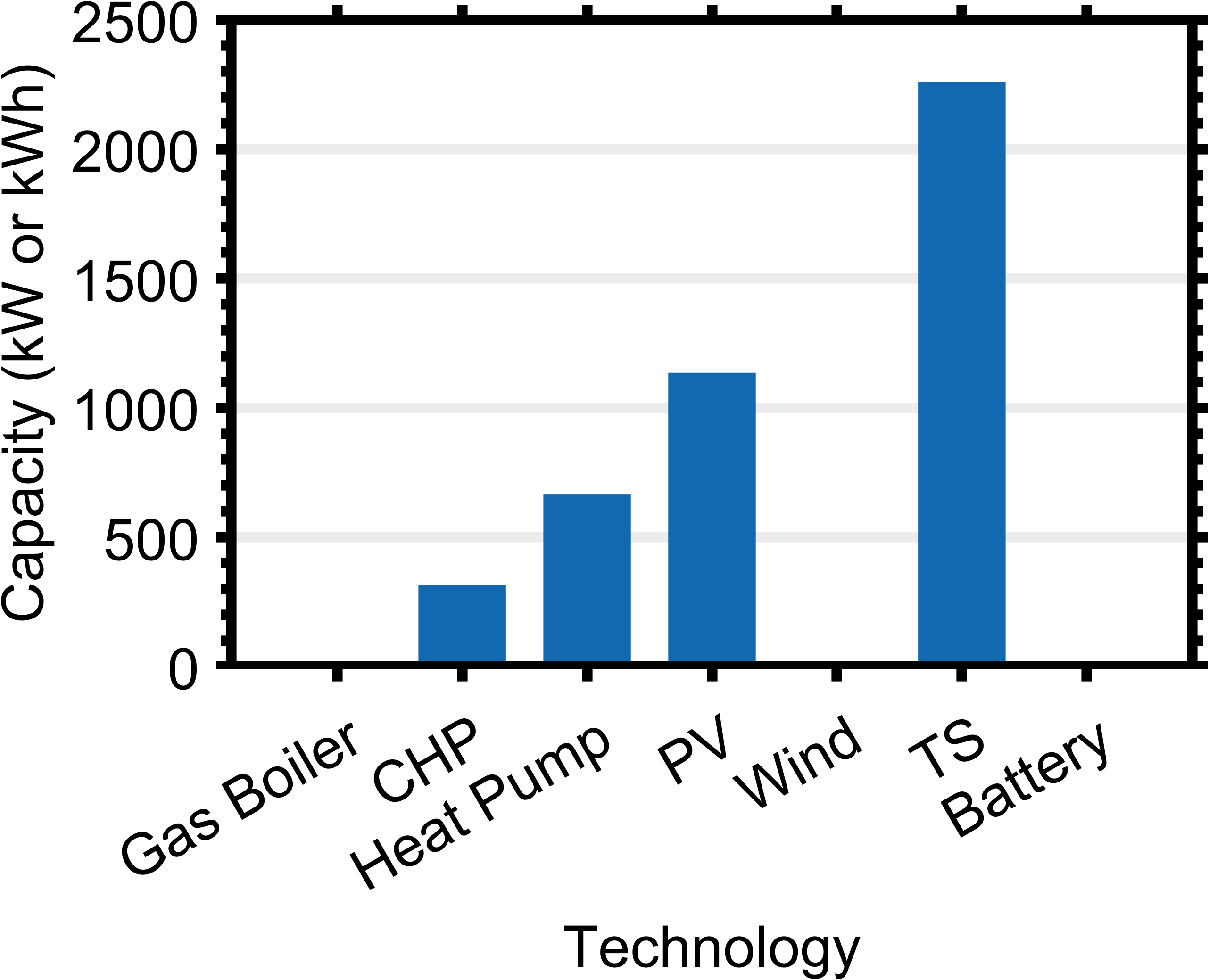}
    \caption{Cost-optimal solution of the Energy Hub model}
    \label{fig:energy_hub_optimal_bar}
\end{figure}

\cref{fig:energy_hub_optimal_bar} shows the solution of the case study, which is optimized for cost subject to the emissions constraint. The combined heat and power unit is the lowest-cost technology for providing heat and electricity. Due to the emissions constraint, however, a substantial amount of \data{PV panels, thermal storage, and a large heat pump} are also included. For further information, see \cref{sec:appendix_cost_optimal}, which shows the operation profile of each technology.

\subsection{Quality of near-optimal spaces}
\label{sec:quality}
\begin{figure}[ht]
    \centering
    \includegraphics[width = \linewidth]{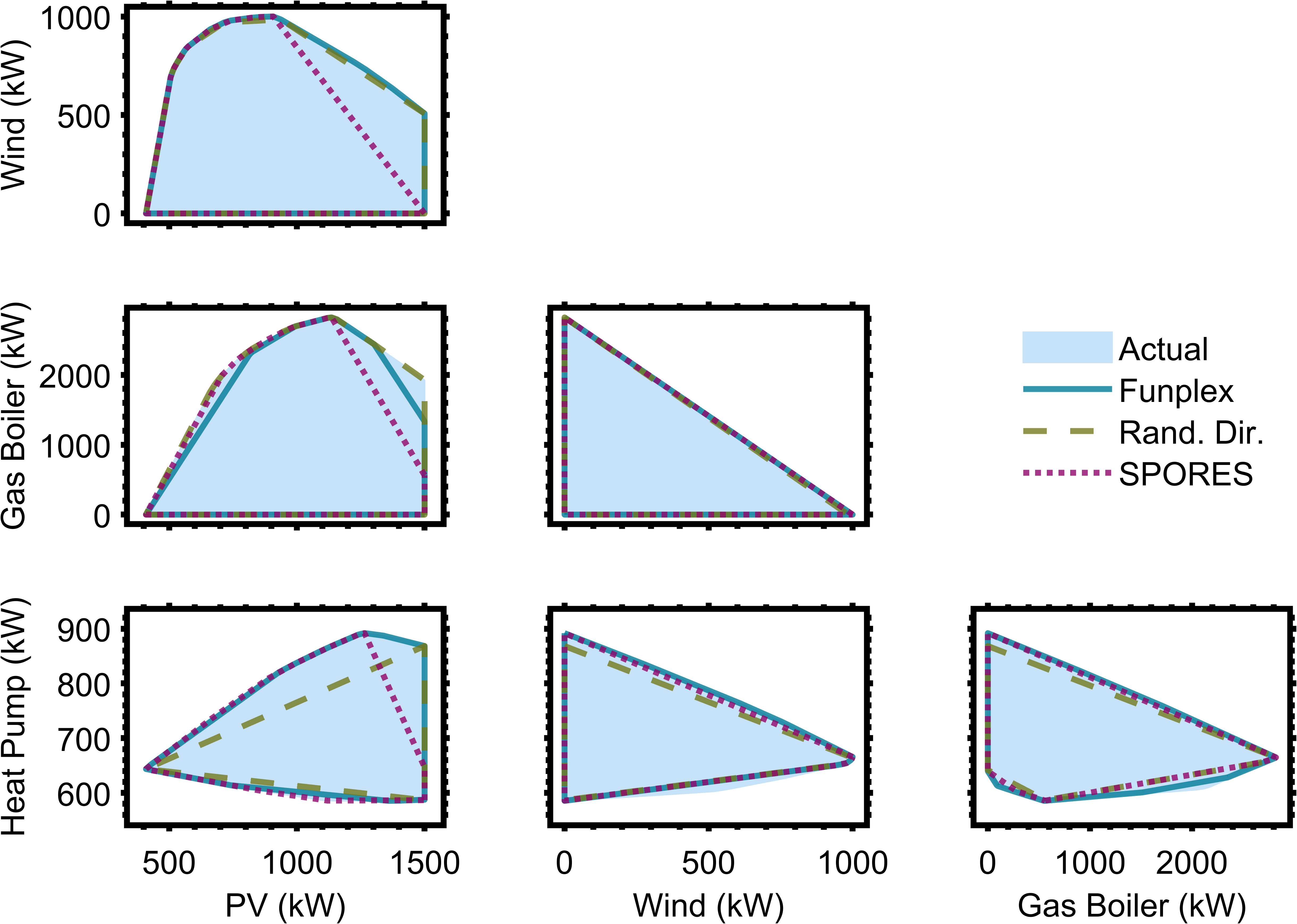}
    \caption{Four-dimensional near-optimal space in terms of wind capacity, PV capacity, gas boiler capacity, and heat pump capacity projected in two dimensions. The lines indicate the outline of the near-optimal space identified by the respective MGA algorithm.}
    \label{fig:MGA_Space_4d}
\end{figure}

\cref{fig:MGA_Space_4d} uses a series of two-dimensional projections to depict the four-dimensional near-optimal space that is found through the MGA analysis. Each panel shows the near-optimal space in terms of two decision variables, which is a customary approach in MGA analyses \cite{grochowicz2023intersecting, neumann2023broad}. The shaded area describes the high-fidelity near-optimal space in these two dimensions obtained via our computationally-intensive planar MGA analysis (see \cref{sec:KPI}), whereas the lines show the outline of the near-optimal space identified using the respective algorithm. An ideal MGA algorithm would delineate the entire blue-shaded region in each panel. 

None of the MGA methods reproduces the space perfectly, accentuating the challenge of MGA to recreate high-dimensional spaces with limited objectives. Nonetheless, Funplex captures most of the near-optimal space correctly and outperforms the other two methods. The region outlined by Funplex closely matches the shaded region in all panels of \cref{fig:MGA_Space_4d}, demonstrating Funplex's effectiveness. SPORES and Random Directions miss larger parts of the near-optimal space for some technologies. The quantitative comparison on \cref{tab:quality} confirms this finding. The volume of the normalized four-dimensional near-optimal space that is found by each algorithm is greatest for Funplex, followed by SPORES. In fact, Funplex offers a volume gain of \data{1.37} and \data{1.25} compared to SPORES and Random Directions, respectively. Hence, Funplex produces higher-quality near-optimal spaces than the established algorithms. 

The improved quality comes from (1) the utilization of intermediary vertices and (2) the choice of objective functions.

\begin{table}[ht]
\centering
\caption{Quality of the near-optimal space found by each of the three algorithms described in terms of the normalized volume of the identified near-optimal space.}
\begin{tabular}{p{0.2\textwidth}>{\centering}p{0.2\textwidth}>{\centering}p{0.2\textwidth}>{\centering\arraybackslash}p{0.2\textwidth}}
\toprule
    & \textbf{Funplex} & \textbf{SPORES} & \textbf{Random Directions} \\ 
\midrule
Volume & \data{1.00} & \data{0.73} & \data{0.80} \\
\bottomrule
\end{tabular}
\label{tab:quality}
\end{table}

\begin{figure}[ht]
    \centering
    \includegraphics[width = \linewidth]{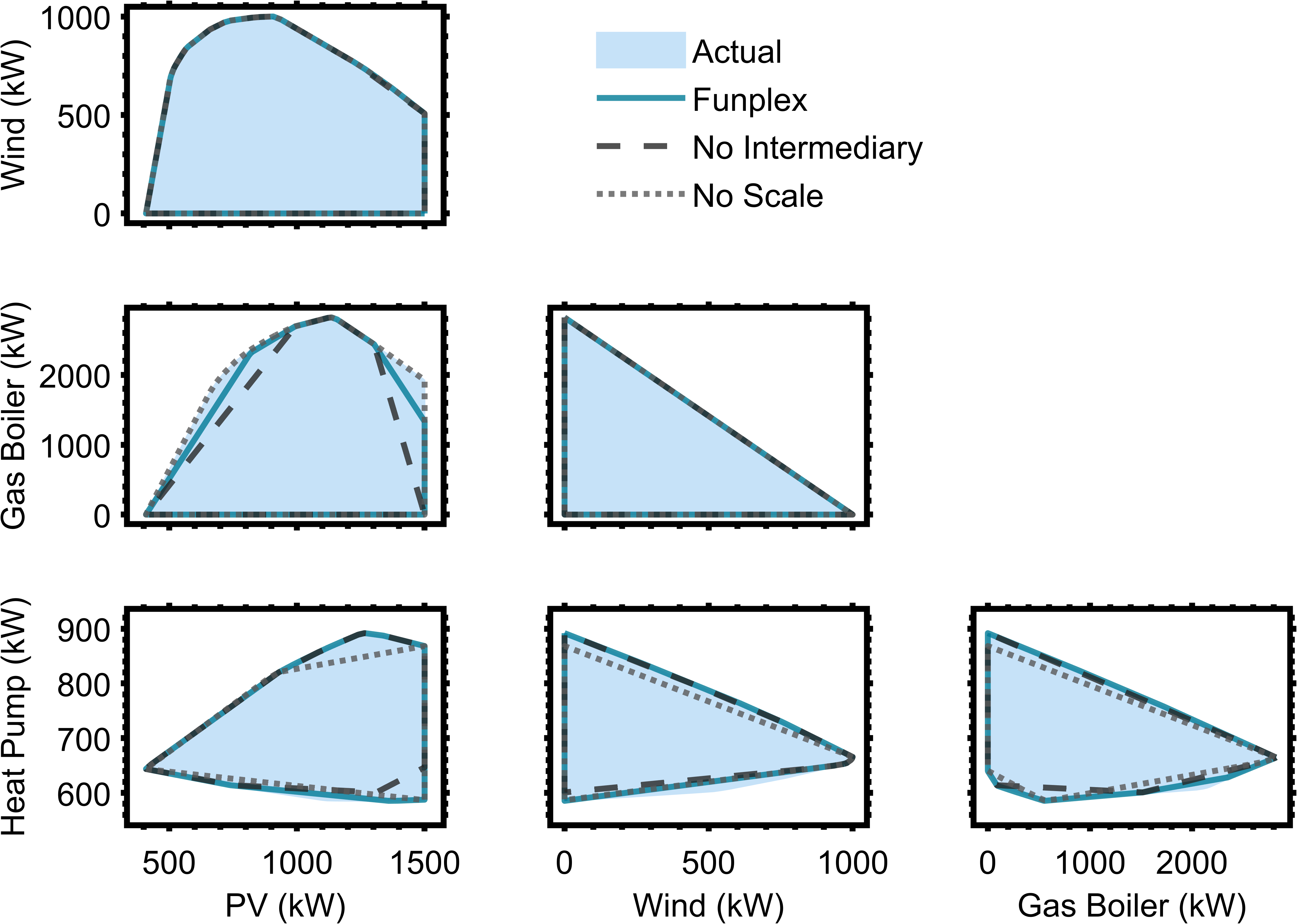}
    \caption{Near-optimal space of the energy-hub model in terms of wind capacity, PV capacity, gas boiler capacity, and heat pump capacity. The plot shows the region identified by Funplex in (1) its standard form, (2) a version which only saves optimal and not intermediary vertices and (3) a version in which variables are not adjusted for scale.}
    \label{fig:efficient_points_only}
\end{figure}

\paragraph{Utilization of intermediary vertices} Funplex saves the coordinates of each intermediary vertex it passes. Doing so helps find larger near-optimal spaces. \cref{fig:efficient_points_only} shows the near-optimal space found by Funplex with (1) both optimal and intermediary vertices and (2) only optimal vertices. The latter region could also be determined by a Random Directions algorithm using the same objective functions as Funplex. \cref{fig:efficient_points_only} illustrates that including intermediary vertices adds information, and a quantitative analysis shows that it increases the volume of the identified near-optimal space by \data{14 \%}. However, this does not explain the full volume gain of Funplex compared to the other two methods. When the intermediary vertices are removed, Funplex still finds higher-quality near-optimal spaces than the other two methods, which is due to a better choice of objective functions.

\clearpage
\paragraph{Choice of objective functions} Funplex's choice of objectives also leads to higher-quality near-optimal spaces. Both Funplex and Random Directions generate the directions randomly; however, objectives differ in two key ways:
\begin{enumerate}
    \item Funplex generates directions from a unit hypersphere. This theoreticaly samples the possible directions more evenly than uniform random numbers. The actual impact is small, increasing the volume of the discovered region by \data{4 \%} in the case study example.\footnote{Some random directions implementations, such as \citeauthor{Patankar2023LandWest} \cite{Patankar2023LandWest}, generate random objective coefficients on the interval between [0,1]. Restricting objective coefficients to be positive in such a way limits the directions with which an algorithm can search the near-optimal space. This greatly impacts the performance of the algorithm, resulting in near-optimal spaces which are over \data{ten} times smaller than the actual spaces on the above test problem. To enable a more fair comparison, we benchmark Funplex against a Random Directions implementation which generates random coefficients on the interval [-1,1] \cite{Berntsen2017EnsuringAlternatives}.}
    \item Funplex accounts for the order of magnitude of the decision variables, allowing it to sample the near-optimal space more evenly. The Random Directions algorithm does not account for the variable scales, thus implicitly favouring variables with larger scales. \cref{fig:efficient_points_only} illustrates that scaling variables increases the volume of the identified near-optimal space by \data{10 \%}. The non-scaled version of Funplex performs better on technologies with large magnitudes such as PV and worse on technologies with smaller scales such as the heat pump. In the case study energy hub model, PV capacity measured in square meters making it larger than other capacity variables. 
\end{enumerate}


\subsection{Computational efficiency}
\label{sec:efficiency}
The indicators benchmarking the three algorithm's performance in terms of computational efficiency are included in \cref{tab:computational_efficiency}. The runtime of Funplex is \data{five} times faster than the other methods when benchmarked using the same underlying Simplex implementation. Thus, Funplex seems to be a promising avenue to improve the computational complexity of MGA and make it more accessible for modeling teams. 

However, whilst runtime is an intuitive indicator of the computational efficiency of an algorithm, it might vary with the exact algorithm implementation. Note that runtime scales linearly with the number of pivots for most implementations. Therefore, the faster runtime comes, in part, due to the lower number of required pivots. Funplex uses around \data{1,900} tableau pivots to determine the space in \cref{fig:MGA_Space_4d}, whereas the two other algorithms require upwards of \data{14,000}. Funplex, therefore, provides an efficiency gain of around \data{eight} compared to both other algorithms. 

Funplex's runtime is hampered slightly by the additional computational cost of checking the optimality of all objectives and saving the current vertex during each pivot. For Simplex, each pivot requires $\mathcal{O}(nm)$ calculations, where $n$ and $m$ are the numbers of constraints and variables, respectively. Since Funplex uses a multi-objective tableau, with additional rows for every objective, the calculations also depend on the number of objectives $N_k$. The complexity of Funplex is  $\mathcal{O}(n(m+N_k))$ per pivot. For small problems like the energy hub case study, the added calculations are significant since $N_k = 200$ is comparable to \data{$m \approx 500$}. The computational penalty is likely negligible for larger problems with $m >> N_k$. The computational complexities depend on the exact simplex implementation. Nonetheless, our analysis suggests that the computational burden of verifying multiple objectives is likely small compared with the burden of each individual pivot. 

Funplex's significant efficiency gain is enabled by (1) warm-starting optimization runs and (2) checking optimality at intermediary vertices.

\begin{table}[ht]
\caption{Computational efficiency of the three algorithms in determining the near-optimal space}
\label{tab:computational_efficiency}
\begin{tabular}{p{0.3\textwidth}>{\centering}p{0.2\textwidth}>{\centering}p{0.2\textwidth}>{\centering\arraybackslash}p{0.2\textwidth}}
\toprule
\textbf{Indicator}    & \textbf{Funplex} & \textbf{SPORES} & \textbf{Random Directions} \\ 
\midrule
Total  Simplex Pivots & \data{1,858}            & \data{19,134}          & \data{14,364}                     \\
Complexity per Pivot  & $\mathcal{O}(n(m+N_k))$       & $\mathcal{O}(nm)$           & $\mathcal{O}(nm)$                      \\
Runtime (s)           & \data{19}               & \data{155}             & \data{125}                        \\ 
\bottomrule
\end{tabular}
\end{table}

\paragraph{Warm-starting optimization runs} Funplex uses information from previously solved objectives to choose efficient starting points for upcoming optimizations. Objectives are optimized in order of similarity. The solver trajectory is reduced since similar objectives often have similar solutions. In contrast, established methods start at the same naive feasible solution regardless of the objective. By only warm-starting the optimization runs without checking intermediary vertices, Funplex requires \data{4117} Simplex pivots. This suggests that warm-starting the optimization runs leads to the largest efficiency improvements compared to the other two algorithms, although checking optimality at intermediary vertices further reduces the number of pivots by about \data{50} \%.

\paragraph{Checking optimality at intermediary vertices} Funplex checks the optimality of all objectives at every vertex it passes. It immediately detects when multiple solutions have the same optimum. Furthermore, it also detects if the optimum of any objective is accidentally passed. Funplex thereby completes some objectives with minimal computational burden. In contrast, other MGA methods always solve complete optimization problems for each objective.

\subsection{Scalability}
\label{sec:scalability}
Funplex is intended as a proof-of-concept, and the results of an initial case study suggest that modified algorithms might improve the quality and computational efficiency of MGA analyses. To understand how Funplex might perform on larger models, \cref{sec:scalability_model} discusses how Funplex's efficiency might scale with model size. Moreover, \cref{sec:scalability_MGA} demonstrates how Funplex scales with MGA complexity to assess its application in more comprehensive analyses.

\subsubsection{Scalability with model size}
\label{sec:scalability_model}
Energy system optimization models typically differ in size through their temporal, spatial, and technological resolution. 

\begin{figure}[h!]
    \centering
    \includegraphics[width=\textwidth]{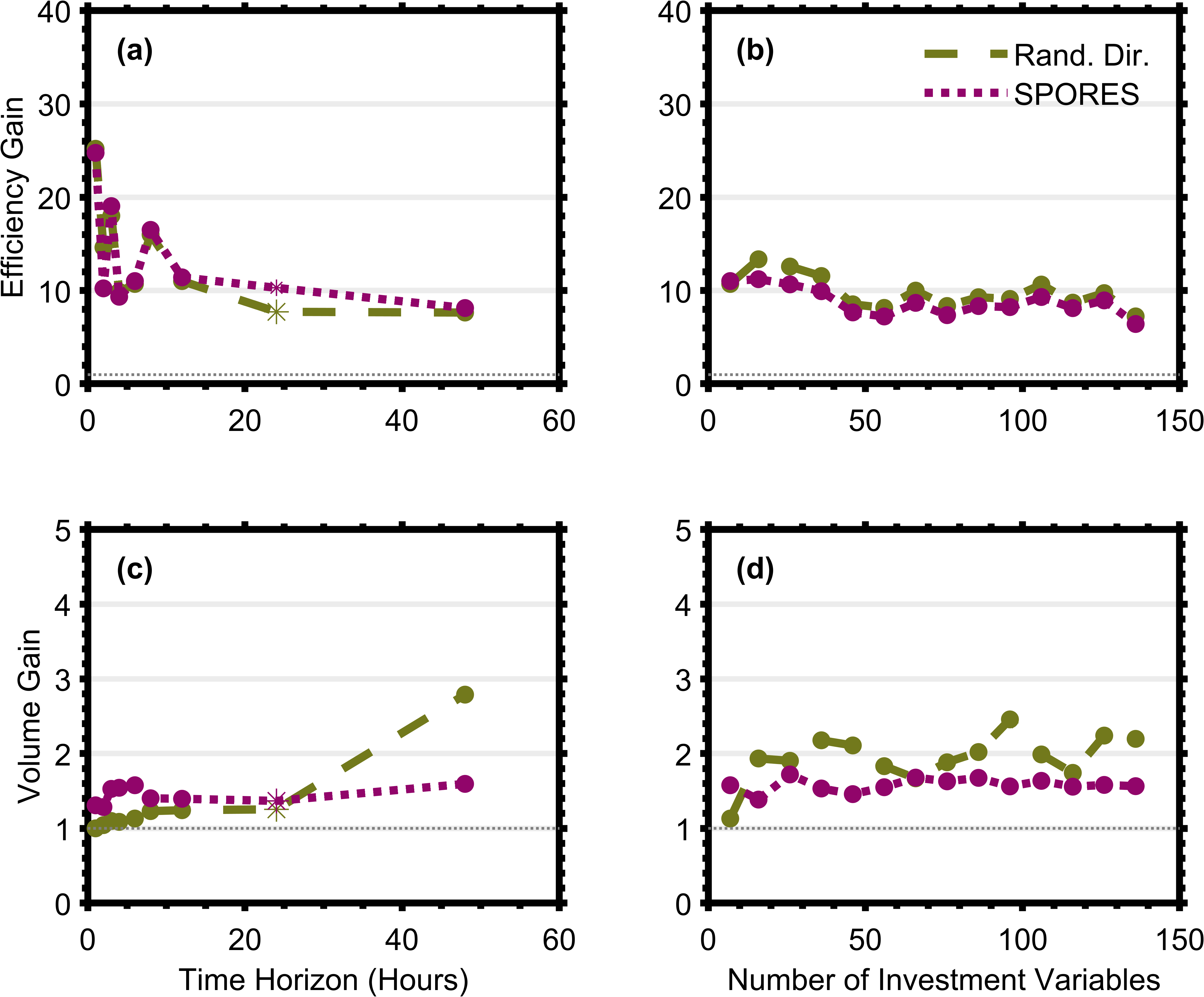}
    \caption{Scalability of Funplex with model size. Panel \textbf{(a)} shows the efficiency gain as a function of the model's operating horizon. Panel \textbf{(b)} shows the efficiency gain as a function of the number of investment variables. Panel \textbf{(c)} shows the volume gain as a function of the model's operating horizon. Finally, Panel \textbf{(d)} shows the volume gain as a function of the number of investment variables. Points marked with an asterisk correspond to the base case shown in \cref{fig:MGA_Space_4d}, \cref{tab:quality}, and \cref{tab:computational_efficiency}. The gray dashed line shows the threshold above which Funplex becomes better than the other algorithms. }.
    \label{fig:runtime_by_problem_size}
\end{figure}

\paragraph{Temporal resolution} A model's temporal resolution determines the number of operational decision variables in a model. Some models simulate only a few representative hours per year \cite{shawhan2021value}, while others simulate all hours in a year. \cref{fig:runtime_by_problem_size}(a) shows how Funplex efficiency changes with the time horizon of the model. As seen in \cref{fig:runtime_by_problem_size}(a), efficiency gain \data{decreases} with the operating horizon. \cref{fig:runtime_by_problem_size}(c) shows that the volume gain compared to Random Directions and SPORES is \data{stable} over the operating horizons tested. There may be a slight \data{upward trend}, which suggests that Funplex achieves greater quality improvements for larger operating horizons. The largest time horizon tested is 48 hours, which is in the same order of magnitude as what is employed in some large, planning-oriented energy-system models \cite{shawhan2021value}. With this time horizon, Funplex still results in an efficiency gain of \data{eight} compared to SPORES and Random Directions. Energy system optimizing models with much larger time horizons, however, are more common. From our scalability analysis, it is unclear whether the efficiency gain plateaus if the time horizon is increased further or whether it continues to diminish. 

\paragraph{Spatial and technological resolution} A model's spatial and technological resolutions determine the number of investment variables in the model. Planning models require separate investment variables for each node and technology simulated. Some models simulate thousands of nodes and many technologies \cite{shawhan2021value}, whereas others simulate only a small number. \cref{fig:runtime_by_problem_size}(b) shows how the efficiency gain changes with the number of investment variables. As seen in \cref{fig:runtime_by_problem_size}(b), the efficiency gain stays \data{approximately constant} with the number of investment variables, which suggests that Funplex may remain fast on problems with many technology options to choose from. \cref{fig:runtime_by_problem_size}(c) shows that the volume gain seems to be \data{independent} of the number of investment variables, confirming that Funplex might be well-suited for models with a large number of investment variables. However, the number of investment variables tested here is lower than typical capacity-planning models, so some uncertainty remains.

\paragraph{Implications of scalability with model size} Funplex might be most suitable for models with a high number of investment variables but a low temporal resolution, such as capacity planning models. Such models include many technologies to choose from but only simulate a small number of representative time periods.  

\subsubsection{Scalability with MGA complexity}
\label{sec:scalability_MGA}
The MGA complexity refers to the resolution and dimensionality of the near-optimal space.

\begin{figure}[h!]
\centering
\includegraphics[width=\textwidth]{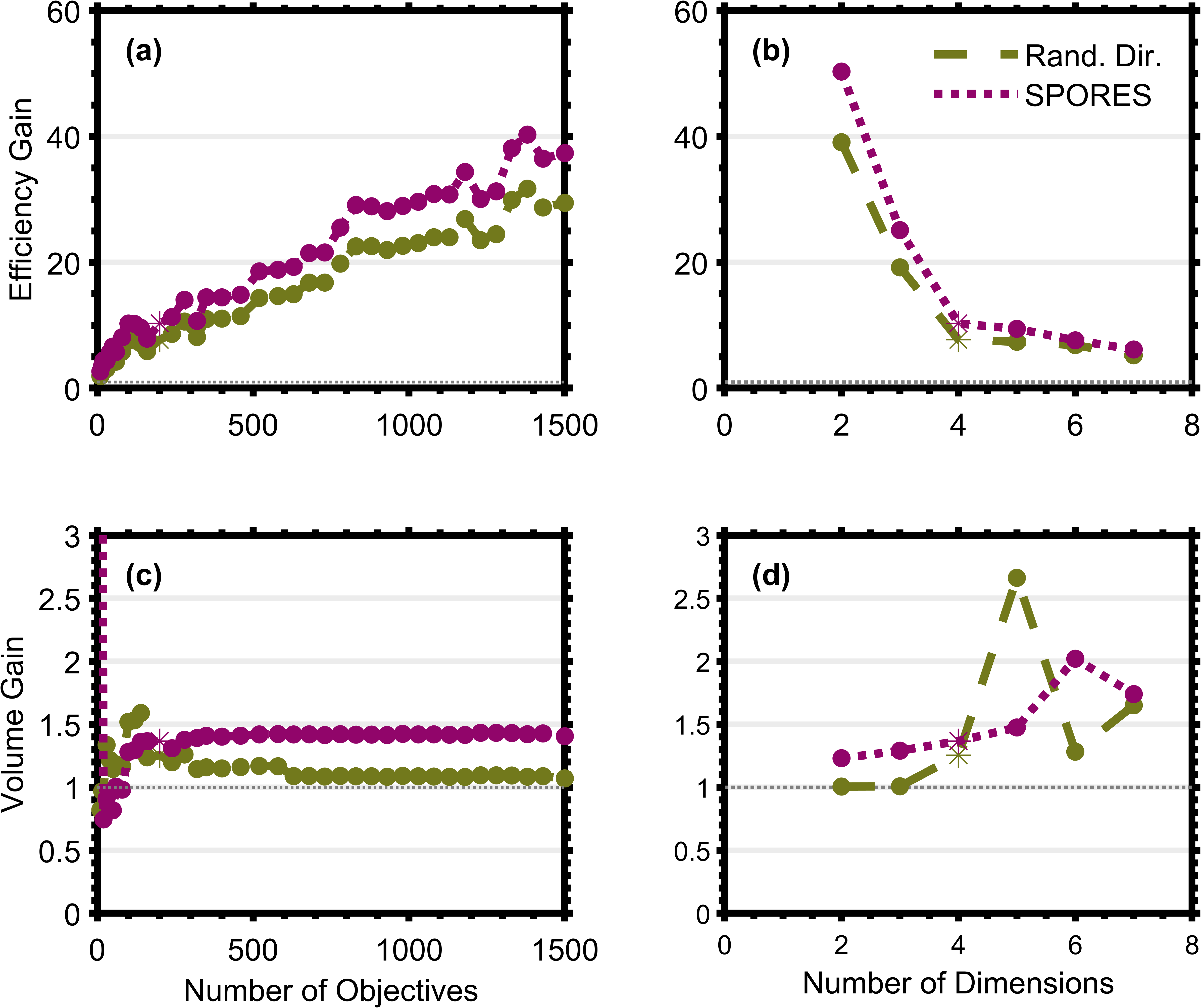}
\caption{Scalability of Funplex with MGA complexity. Panel \textbf{(a)} shows the efficiency gain as a function of MGA objectives. Panel \textbf{(b)} shows the efficiency gain as a function of the dimensions of the MGA analysis. Panel \textbf{(c)} shows the volume gain as a function of MGA objectives. Panel \textbf{(d)} shows the volume gain as a function of the dimension of the MGA analysis. Points marked with an asterisk correspond to the base case shown in \cref{fig:MGA_Space_4d}, \cref{tab:quality}, and \cref{tab:computational_efficiency}. The gray dashed line shows the threshold above which funplex becomes better than the other algorithms. }
\label{fig:runtime_by_complexity}
\end{figure}

\paragraph{MGA objectives} Near-optimal spaces are of higher quality when the determined area closely matches the high-fidelity near-optimal space. The resolution is improved by optimizing more objectives. \cref{fig:runtime_by_complexity}(a) shows that Funplex efficiency \data{increases} strongly with the number of optimized objectives. In Funplex, the number of pivots scales roughly \data{$\mathcal{O}(\sqrt{N_k})$},\footnote{Determined through a log-log regression} whereas both SPORES and Random directions scale linearly.  Next, as seen in \cref{fig:runtime_by_complexity} (c) the volume gain compared to SPORES remains \data{constant} for increasing numbers of MGA objectives. The volume gain of Random Directions converges towards one, since both Funplex and Random Directions find the entire near-optimal space for a large number of objectives.   Thus, overall, Funplex is increasingly more efficient with more objectives optimized. With more objectives, conventional methods will likely perform more optimizations that repeat boundary points. Funplex eliminates this inefficiency by checking the optimality of all objectives at every vertex. Similarly, Funplex is more likely to start close to the desired optimum. Funplex is, therefore, good at finding near-optimal spaces with a high resolution.

\paragraph{MGA dimensions} Near-optimal spaces are high-dimensional when many decision variables of interest are considered simultaneously. High-dimensional spaces are difficult to determine because the number of variable combinations increases rapidly with the number of dimensions \cite{grochowicz2023intersecting}. Algorithms such as SPORES are designed to work with high-dimensional spaces efficiently. \cref{fig:runtime_by_complexity}(b) shows how efficiency changes with the dimensions of the near-optimal space. As seen, efficiency drops with dimension. Spaces in higher dimensions are larger and more diverse. Objectives are, therefore, less likely to find the same boundary point, and the difference between optimal solutions of various objectives is greater. Consequentially, there are fewer computational redundancies which Funplex can exploit in larger dimensions. However, \cref{fig:runtime_by_complexity}(d) suggests that whilst the efficiency gain might decrease for an increasing number of dimensions, the quality of the found solutions \data{increases} slightly. Thus, Funplex might still be a valuable option for high-dimensional analyses if a high-quality solution is desired. The largest number of dimensions tested was seven because this is the maximum number of technologies in the energy hub model. Seven dimensions are the same order of magnitude as in many MGA analyses \cite{grochowicz2023intersecting, neumann2023broad}, and Funplex is still roughly \data{five} times more efficient at this level whilst finding higher-quality solutions.

\paragraph{Implications of scalability with MGA complexity}
Funplex might provide the largest computational efficiency gain for studies with low-dimensional near-optimal spaces, such as the studies by Grochowicz et al. \cite{grochowicz2023intersecting} and Neumann and Brown \cite{neumann2023broad}. Both studies identify all possible technology mixes for a future European energy system. The near-optimal spaces they analyze only have five to ten dimensions total, each representing the total installed capacity of a given key technology. Higher-dimensional analyses, such as the one performed by Lombardi et al. \cite{lombardi2020policy} assessing the spatial distribution of many technologies, might yield less computational efficiency gains, although the quality of the found solutions might still be increased. 
Our analysis suggests that overall, Funplex is a promising method when it is important to capture all, or nearly all, of the near-optimal space since Funplex produces higher-quality near-optimal spaces than current methods and scales well with the number of objectives optimized, especially compared to established MGA algorithms. Funplex may, therefore, have similar applications to Modeling All Alternatives, an algorithm that seeks to determine the complete near-optimal space \cite{pedersen2021modeling} whilst being potentially more scalable when built on more advanced solvers. Since the additional computational burden of optimizing more objectives is small, Funplex reduces the urgency of picking objectives wisely. Users can optimize more objectives rather than fine-tune MGA methods and parameters, making the algorithm user-friendly compared to existing methods. 

\subsection{Limitations}
\label{sec:limitations}
The current implementation of Funplex is intended as a proof of concept. It is numerically unstable for problems with more than \data{1000} decision variables, thus limiting the extent to which a scalability analysis can be performed. If the numerical stability is improved, memory limitations might become a limiting factor to significantly increasing the problem size since the entire tableau is stored at each iteration. By building on more advanced algorithms, such as the revised Simplex or dual Simplex algorithm, such problems might be overcome. These algorithms start every iteration with the original coefficients, thus improving numerical stability. Moreover, they use sparse matrices to store the tableau, thus reducing the memory requirements. We suggest advancing the algorithm in collaboration with professional solver developers to overcome these limitations and to further improve the algorithm speed. 

\section{Conclusions}
\label{sec:Conclusions}
As an increasingly popular method in energy systems optimization, modeling to generate alternatives (MGA) enables the exploration of near-cost-optimal system alternatives that stakeholders might prefer for various reasons. Current MGA methods treat solvers as black boxes, and efforts to improve the computational efficiency of MGA focus on strategically choosing objective functions. Each objective is optimized separately, resulting in repeated calculations and computational inefficiencies. There has been little to no research on whether modified solver algorithms can reduce these inefficiencies, and there is a poor understanding of the achievable performance improvements. This paper attempts to address these research gaps via the proof-of-concept algorithm Funplex. 

Funplex is a Simplex-based algorithm specifically designed to find near-optimal spaces efficiently. The algorithm generates random MGA objectives from a unit hypersphere and uses a multi-objective Simplex tableau to optimize the objectives efficiently. Funplex considers information from all objectives throughout the whole solving process, thereby minimizing computational redundancies. 

Applied to a case study energy hub model, Funplex produces higher-quality near-optimal spaces and is \data{five} times faster than established MGA methods. It produces higher-quality spaces by (1) generating better MGA objectives and (2) saving the coordinates of intermediary vertices along its trajectory. Funplex is faster than established methods because it (1) checks the optimality of all objectives at each passed vertex and (2) uses information from previously solved optimizations to initialize future optimizations efficiently.

A scalability analysis suggests that Funplex scales well with investment variables but shows a diminishing performance for increasing operating horizons. Thus, Funplex might be most suitable for large capacity planning problems. Moreover, Funplex's efficiency gains increase with the number of MGA objectives but decrease with the number of MGA dimensions, making Funplex most promising for lower-dimensional analyses requiring good coverage of the near-optimal space.

Overall, our study provides strong evidence that new solver algorithms may significantly improve the quality and computational efficiency of MGA. This analysis is intended as a proof-of-concept. Funplex could be further improved by collaborating with professional solver developers to build on more advanced algorithms such as the dual Simplex, revised Simplex, or other non-Simplex-based algorithms. Such improvements would allow for testing the algorithm on larger models, helping to determine further aspects for improvement.

\section*{Acknowledgement}
This project has received funding from the European Union’s Horizon 2020 research and innovation program under the Marie Sklodowska-Curie (MSC) grant agreement No.847585 - RESPONSE.


\appendix


\section{Simplex algorithm}
\label{sec:appendix_Simplex}

Simplex is a standard method for solving linear optimization problems. The algorithm moves along the boundary of the feasible region, from vertex to vertex. It does so in a way that always moves closer to the optimal point. Since the number of vertices is finite, the optimum is eventually reached. 

Consider a single-objective linear optimization problem in standard optimization form:  

\begin{equation}
\label{eq:standard_form}
\begin{aligned}
\min_{\mathbf{x}} \quad & \mathbf{c}^T \mathbf{x}\\
\textrm{s.t.} \quad & A \mathbf{x} = \mathbf{b}\\
  &\mathbf{x} \geq \mathbf{0}    \\
\end{aligned}
\end{equation}

\noindent Let $n$ be the number of decision variables in the problem and $m$ be the number of constraints. Then, $\mathbf{x} \in \mathbb{R}^n$ is a vector of decision variables; $\mathbf{c} \in \mathbb{R}^n$ is a vector of objective coefficients; $A \in \mathbb{R}^{m,n}$, is a matrix of constraint coefficients which is assumed to be of full rank; and $\mathbf{b} \in \mathbb{R}^{m,1}$ is a vector of constraint constants.

A core concept in Simplex is the notation of a \textit{basic feasible point} \cite{Nocedal2006NumericalOptimization}. Let $\mathcal{F}$ denote the feasible region of the optimization problem \cref{eq:standard_form}. Then, a basic feasible point is a point $\tilde{\mathbf{x}} \in \mathcal{F}$ such that there exists a set of indices $\mathcal{B} \subset \{1, ..., n\}$ where:

\begin{enumerate}
    \item $\mathcal{B}$ contains exactly $m$ indices.
    \item $\tilde{x}_i = 0$ for every $i \not \in \mathcal{B}$.
    \item The columns of $A$ whose indices are given by $\mathcal{B}$ are linearly independent.
    
\end{enumerate}

\noindent The set $\mathcal{B}$ is said to be a basis of \cref{eq:standard_form}. For further details, see \citeauthor{Nocedal2006NumericalOptimization} \cite{Nocedal2006NumericalOptimization}.

Feasible basic points have a few key properties which are exploited by Simplex. First, every vertex of $\mathcal{F}$ is also a feasible basic point. Next, if \cref{eq:standard_form} has a feasible point (i.e. $\mathcal{F} \neq \emptyset$), then it also has a feasible basic point. Finally, if \cref{eq:standard_form} has a solution, then at least one of solution is a feasible basic point. Proofs of all three statements are provided in \citeauthor{Nocedal2006NumericalOptimization} \cite{Nocedal2006NumericalOptimization}. These properties imply that the optimal solution can be found by searching only through the set of feasible basic points. Solving \cref{eq:standard_form} is thus reduced to a new problem of finding the optimal basis. 

Simplex systematically solves this problem by partitioning the decision variables into basic ($\mathbf{x}_B \in \mathbb{R}^m$) and non-basic ($\mathbf{x}_N \in \mathbb{R}^{n-m}$) variables. Basic variables can take on any value, while non-basic variables are all zero ($\mathbf{x}_N = \mathbf{0}$). The basic variables should be chosen to that their indices form the basis of a fesaible basic point. When this is the case, the partition defines a vertex of the feasible region.  After separating the decision variables, the problem formulation (\cref{eq:standard_form}) can be decomposed as \cite{Hug2023ETHSystems}: 

\begin{equation}
\label{eq:decomposed_standard_form}
\begin{aligned}
\min_{\mathbf{x}_B, \mathbf{x}_N} \quad & \mathbf{c}_B^T \mathbf{x}_B + \mathbf{c}_N^T \mathbf{x}_N\\
\textrm{s.t.} \quad & B \mathbf{x}_B + D \mathbf{x}_N = \mathbf{b}\\
  &\mathbf{x}_B, \mathbf{x}_N\geq \mathbf{0}    \\
\end{aligned}
\end{equation}

\noindent where $B \in \mathbb{R}^{m,m}$, $D \in \mathbb{R}^{m,n-m}$, $\mathbf{c}_B \in \mathbb{R}^{m}$, and $\mathbf{c}_N \in \mathbb{R}^{n-m}$. 

Simplex uses this decomposed problem (\cref{eq:decomposed_standard_form}) to define the Simplex tableau. The tableau contains all the information of the original problem but makes it easier to read important problem characteristics. For a given basis, the Simplex tableau is defined as \cite{Hug2023ETHSystems}:

\begin{equation}
\label{eq:Simplex_tableau}
   T = \left[\arraycolsep=5pt \def\arraystretch{1.2}
    \begin{array}{cc|c} 
	I & B^{-1} D & B^{-1} \mathbf{b}\\ \hline 
	0 & \quad \mathbf{c}_{D}^{T} - \mathbf{c}_{B}^{T} B^{-1} D & -\mathbf{c}_{B}^{T} B^{-1} \mathbf{b}
\end{array}\right] 
\end{equation}

The top half of the tableau, with dimensions  $m \times (n+1)$, contains the constraints of the original problem. The constraints have been manipulated via elementary row operations to make the current vertex coordinates easy to read. The tableau describes a vertex of the feasible region. At this vertex, the values of the basic variables are given by the top-right quadrant of the tableau. The full vertex coordinates are $\mathbf{x}_N = 0$ and $\mathbf{x}_B = B^{-1}b$.

In contrast, the bottom half of the tableau contains the objective function, with dimensions $1 \times (n+1)$. The bottom right quadrant is the negative of the objective value at the current vertex. The bottom-left quadrant is called the \emph{relative-cost vector}. This vector describes how the objective function will change when each non-basic variable is increased by one unit. 

Once an initial tableau has been created, Simplex uses the relative cost vector to determine what vertex to move to next. For any non-basic variable with a negative relative cost, the objective will decrease when the variable is made non-zero. Simplex can thus improve the objective value by switching this non-basic variable with a basic variable \cite{zenklusen2023linear}. Simplex switches chosen non-basic and basic variables by \emph{pivoting} the Simplex tableau. This involves using elementary row operations to put the tableau into the proper form for the new basis \cite{Hug2023ETHSystems}. Simplex repeats this pivoting process many times to improve the objective. Importantly, the optimal solution has been obtained when all elements of the relative cost vector are either zero or positive. At this point, the objective value can no longer be improved \cite{zenklusen2023linear}. 

\section{Established MGA methods}
\label{sec:appendix_MGA_methods}
Two established MGA methods are relevant for this paper: SPORES \cite{lombardi2020policy, lombardi2023redundant} and Random Directions \cite{patankar2023land}. These two methods differ in the variables of interest ($\mathbf{x^{*}}$) they consider and the objective weights ($\mathbf{w}$) they use. The SPORES method is iterative, meaning that the objective weights of the next optimization depend on the solution of the previous optimizations. In contrast, the Random Directions method is non-iterative: all objectives can be determined immediately at the start of the algorithm.

\subsection{SPORES}
The SPORES method explores spatially-resolved near-optimal spaces \cite{lombardi2020policy, lombardi2023redundant}. Let $\mathcal{I}$ be the set of model technologies and  $\mathcal{J}$ be the set of model regions. The variable of interest in SPORES is the installed capacity of technology $i$ in a given region $j$, i.e., $\mathbf{x^{*}} = \{x_{i,j}^{cap} | \:i \in \mathcal{I}, \:j\in \mathcal{J}\}$. Spatially-resolved near-optimal spaces are high-dimensional and, therefore, difficult to explore thoroughly. In each iteration, SPORES attempts to gain as much information as possible by choosing an objective that generates a maximally different solution from the previous ones. The SPORES MGA formulation is: 

\begin{equation}
\label{eq:SPORES}
\begin{aligned}
\min_{\mathbf{x}} \quad & a \sum_{j\in J} (x_{i_0 j}^{cap}) + b\sum_j \sum_i w_{ij} (x_{ij}^{cap})\\
\textrm{s.t.} \quad & A \mathbf{x} = \mathbf{b}\\
  & \mathbf{c}^T \mathbf{x} \leq (1+\epsilon)f_{min} \\ 
  &\mathbf{x} \geq \mathbf{0}    \\
\end{aligned}
\end{equation}

The objective function has two different terms, weighted by parameters, $a$ and $b$. The first term is used to favor or hamper a single technology of interest ($i_0 \in \mathcal{I}$). It assigns a weight to the total aggregate capacity of that technology over all regions. This guarantees solutions which technologically diverse, and ensures that all technological options are fully explored. The second term incentivizes solutions that are technologically and spatially different from previous ones. 

The variable weights $w_{ij}$ are determined iteratively using the following recursion:

\begin{equation}
\label{eq:spores_innitialization}
    w_{ij}^0 =     \begin{cases}
        r_0, \quad\textup{if}\:\: x_{ij}^{cap,opt} > \gamma\\
        0, \:\:\quad\textup{if}\:\: x_{ij}^{cap,opt} \leq \gamma\\
    \end{cases}
\end{equation}

\begin{equation}
\label{eq:spores_recursion}
    w_{ij}^n = w_{ij}^{n-1} + r_{ij}, \quad \textup{where}\:\:
    r_{ij} =
    \begin{cases}
        r_0, \quad\textup{if}\:\: x_{ij}^{cap, n-1} > \gamma\\
        0, \:\:\quad\textup{if}\:\: x_{ij}^{cap, n-1} \leq \gamma\\
    \end{cases}
\end{equation}

Here, $x_{ij}^{cap, opt}$ and $x_{ij}^{cap, n}$ are the values of $x_{ij}^{cap}$ in the cost-optimal solution and $n^{th}$ iteration, respectively. Next, $w_{ij}^n$ is the weight $w_{ij}$ in the $n^{th}$ objective. The parameters $r_0$ and $\gamma$ can be tuned to match the problem formulation and are chosen to be \data{$r_0 = 0.5$} and \data{$\gamma = 100$} kW for this analysis.  The recursion in \cref{eq:spores_recursion,eq:spores_innitialization} is one of several options provided in \cite{lombardi2023redundant} and is chosen here by recommendation of the authors.

SPORES runs several optimization sequences in parallel. Each sequence is defined by a technology $i_0$, and a set of coefficients ($a$, $b$). Each individual sequence rapidly converges to a fixed point, where $x_{ij}^{n-1} = x_{ij}^{n}$. Using multiple sequences therefore allows the user to continue searching the near-optimal space. In this analysis, we use the parameter values \data{$(a, b) \in \{(-100, 1), (-10, 1), (-1, 1), (1, 1), (10, 1), (100, 1)\}$} These values are similar to the ones chosen in \cite{lombardi2020policy, lombardi2023redundant}. We weight the first objective term more heavily since the analyses in this paper focus on technological near-optimal spaces and not spatial ones. The total number of simulations run by SPORES is therefore \data{$6 \times \lvert\mathcal{I}\rvert \times N_{max}$}, where $\lvert\mathcal{I}\rvert$ is the number of technologies and $N_{max}$ is the number simulations in each sequence. The coefficient \data{6} comes from the number of parameter combinations ($a$, $b$) simulated.  

\subsection{Random Directions algorithm}
The Random Directions method identifies the near-optimal space using randomly generated objectives. There are multiple implementations, but this paper relies on the recently published study by \citeauthor{Berntsen2017EnsuringAlternatives} \cite{Berntsen2017EnsuringAlternatives}. The variables used in the MGA formulation are $\mathbf{x^{*}} = \{x_{i}^{*} | \:i \in \mathcal{I}\}$ where $x_{i}^{*}$ is the $i^{th}$ variable of interest.

The objective weights of these variables aredetermined randomly on the interval [-1,1]. Let $\beta_i \sim \textup{Unif}(-1,1)$ be such a randomly generated objective coefficient. The Random Directions MGA formulation is then:

\begin{equation}
\label{eq:random_directions_formulation}
\begin{aligned}
\max/\min_{\mathbf{x}} \quad & \sum_i \beta_i x_{i}^{prod}\\
\textrm{s.t.} \quad & A \mathbf{x} = \mathbf{b}\\
  & \mathbf{c}^T \mathbf{x} \leq (1+\epsilon)f_{min} \\ 
  &\mathbf{x} \geq \mathbf{0}    \\
\end{aligned}
\end{equation}

Random directions repeatedly solves \cref{eq:random_directions_formulation} to obtain different boundary points. Unlike SPORES, it is not iterative. Each optimization problem is completely independent of the previous ones.

\section{Cost optimal solution}
\label{sec:appendix_cost_optimal}

\begin{figure}[h!]
    \centering
    \includegraphics[width=\textwidth]{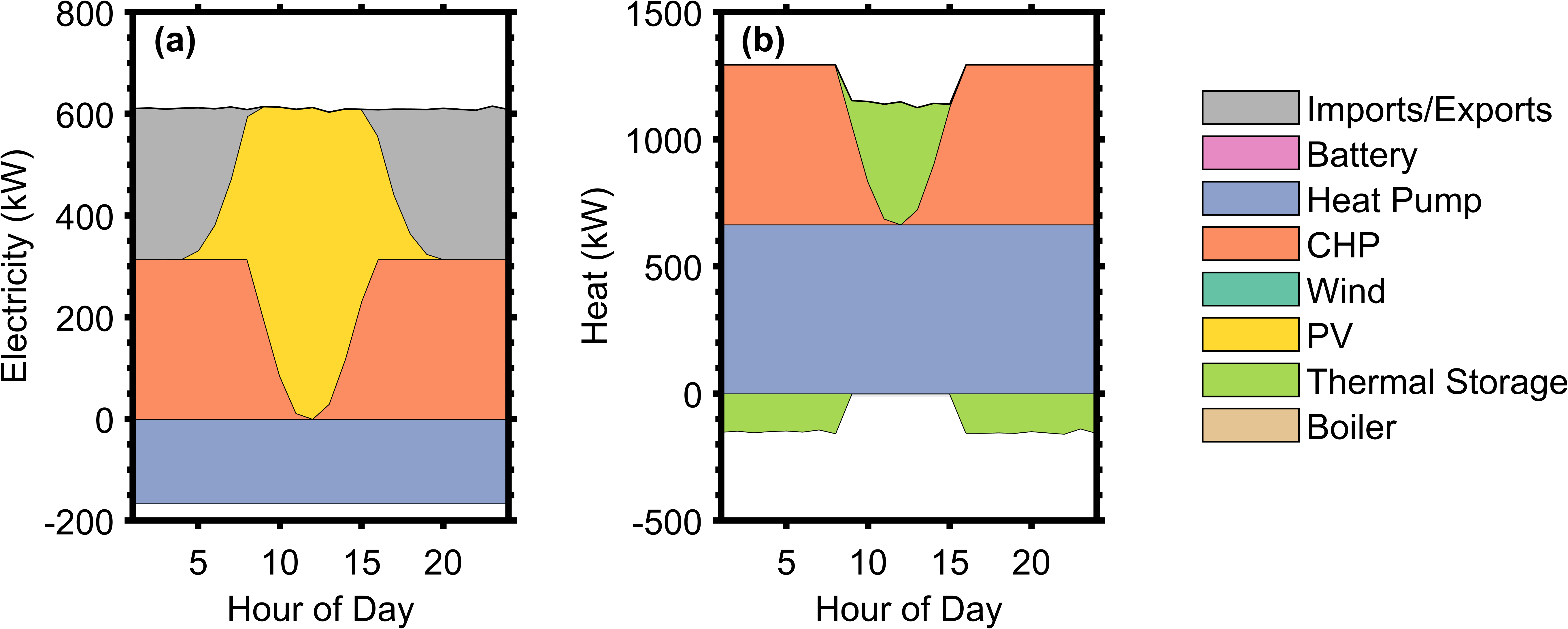}
    \caption{Operation profile of technologies in the energy hub model. Panel \textbf{(a)} shows electricity produced and consumed by each generating technology throughout the representative day, and panel \textbf{(b)} shows the heat produced and consumed by each heating technology throughout the representative day.}
    \label{fig:stackplot}
\end{figure}

\cref{fig:stackplot} shows the operation profile of the cost-optimal solution in the energy hub model. Panel (a) shows the operation profile of each electricity-producing technology. During the night, electricity is produced by the CHP unit and supplemented through grid imports. Then, during the day, solar generation replaces first imports and then CHP. The electricity is consumed by a heat pump, which operates in a baseload fashion. The remaining electricity is used to meet the electricity demand. Panel (b) shows the operation profile of all heating technologies. During the night, heat demand is met through output of the CHP unit and the heat pump. During the day, TS is used to supply the heating demand while the CHP unit is operating below capacity. 

\newpage
\bibliographystyle{elsarticle-num-names} 
\bibliography{references}





\end{document}